\documentstyle{amsppt}
\voffset-10mm
\magnification1200
\pagewidth{130mm}
\pageheight{208mm}
\hfuzz=2.5pt\rightskip=0pt plus1pt
\binoppenalty=10000\relpenalty=10000\relax
\TagsOnRight
\loadbold
\nologo
\addto\tenpoint{\normalbaselineskip=1.05\normalbaselineskip\normalbaselines}
\addto\eightpoint{\normalbaselineskip=1.05\normalbaselineskip\normalbaselines}
%==================================================
\let\hlop\!

\redefine\d{\roman d}
\redefine\Re{\operatorname{Re}}
\redefine\pmod#1{\;(\operatorname{mod}#1)}
\define\SL{SL}
\define\OO{O}
%==================================================
\topmatter
\title
Well-poised generation \\
of Ap\'ery-like recursions
\endtitle
\author
Wadim Zudilin\footnotemark"$^\ddag$"\ \rm(Moscow)
\endauthor
\date
\hbox to70mm{\vbox{\hsize=70mm%
\centerline{E-print \tt math.NT/0307058}
\smallskip
\centerline{May/November 2003}
}}
\enddate
\address
\hbox to70mm{\vbox{\hsize=70mm%
\leftline{Moscow Lomonosov State University}
\leftline{Department of Mechanics and Mathematics}
\leftline{Vorobiovy Gory, GSP-2}
\leftline{119992 Moscow, RUSSIA}
\leftline{{\it URL\/}: \tt http://wain.mi.ras.ru/index.html}
}}
\endaddress
\email
{\tt wadim\@ips.ras.ru}
\endemail
\abstract
The idea to use classical hypergeometric series and, in particular,
well-poised hypergeometric series in diophantine problems
of the values of the polylogarithms has led to several novelties
in number theory and neighbouring areas of mathematics.
Here we present a systematic approach to derive
second-order polynomial recursions for approximations
to some values of the Lerch zeta function, depending on the fixed
(but not necessarily real) parameter $\alpha$
satisfying the condition $\Re(\alpha)<1$.
Substituting $\alpha=0$ into the resulting recurrence
equations produces the famous recursions for
rational approximations to $\zeta(2)$, $\zeta(3)$ due to Ap\'ery,
as well as the known recursion for rational approximations to~$\zeta(4)$.
Multiple integral representations for solutions
of the constructed recurrences are also given.
\endabstract
% \subjclass
% 11J70, 33C20, 33F10
% \endsubjclass
\endtopmatter
%\footnote""{2000 {\it Mathematics Subject Classification}.\enspace
%11J70, 33C20, 33F10.}
\leftheadtext{W.~Zudilin}
\rightheadtext{Well-poised generation of Ap\'ery-like recursions}
\footnotetext"$^\ddag$"{The work is supported by an Alexander von Humboldt
research fellowship.}
%==================================================
\document

The idea to use classical hypergeometric series \cite{9, 11}
and, in particular, well-poised hypergeometric series \cite{13}
in diophantine problems of the values of the polylogarithms has led
to several novelties in number theory and
neighbouring fields of mathematics.
Here we present a systematic approach to derive second-order
polynomial recursions for approximations to the numbers
$$
Z_2^-(\alpha)=\sum_{\nu=1}^\infty\frac{(-1)^{\nu-1}}{(\nu-\alpha)^2},
\quad
Z_3(\alpha)=\sum_{\nu=1}^\infty\frac1{(\nu-\alpha)^3},
\quad
Z_4(\alpha)=\sum_{\nu=1}^\infty\frac1{(\nu-\alpha)^4},
$$
where the fixed (but not necessarily real) parameter~$\alpha$
satisfies the condition $\Re(\alpha)\hlop<\hlop1$.
Substituting $\alpha=0$ into the resulting recurrence equations
produces the famous recursions for rational approximations to
$\zeta(2)=2Z_2^-(0)$, $\zeta(3)=Z_3(0)$ due to Ap\'ery~\cite{1},
as well as the recursion for rational approximations
to $\zeta(4)=Z_4(0)$ known as the Cohen--Rhin--Sorokin--Zudilin
recursion (\cite{8, 15, 17}), which is proper
from both the historical and the alphabetic point of view.

To make clear to the reader, what do we mean by a {\it recursion
for approximations to a number\/} $z\in\Bbb C$, we introduce a
formal definition. The requirement to such the recursion is to
have two linearly independent solutions $\{u_n\}_{n=0}^\infty$
and $\{v_n\}_{n=0}^\infty$ (uniquely determined by
the recursion itself and initial conditions) such that
$v_n/u_n\to z$ as $n\to\infty$. In the case $z\in\Bbb R$
(e.g., corresponding to $\alpha\in\Bbb Q$ in the above definitions),
we usually restrict ourselves to the sequences
$\{u_n\}_{n=0}^\infty$, $\{v_n\}_{n=0}^\infty$
consisting of rational numbers only and may interpret
them as denominators and numerators, respectively, of
rational convergents to~$z$.

We apologise in advance for facing the reader with,
sometimes, cumbersome formulae. Although the ideas
of the well-poised hypergeometric construction
of linear forms $r_n=u_nz-v_n$, $n=0,1,2,\dots$, are simple,
the appearance of lengthy formulae is unavoidable
in this type of analysis. Our main results
are the recursions \thetag{4}, \thetag{8},
\thetag{10}, \thetag{12}, \thetag{13},
and the integral representations \thetag{6},
\thetag{9}, \thetag{11} (generalizing those of~\cite{3})
for the linear forms $\{r_n\}_{n=0}^\infty$.

\subhead
1. $Z_2^-(\alpha)$
\endsubhead
Take the rational function
$$
\align
R_n(t)
&=n!\cdot(2t-2\alpha+n)
\frac{(t-1)\dotsb(t-n)\cdot(t-2\alpha+n+1)\dotsb(t-2\alpha+2n)}
{\bigl((t-\alpha)(t-\alpha+1)\dotsb(t-\alpha+n)\bigr)^3}
\\
&=2n!\cdot\frac{\Gamma(t-\alpha+n/2+1)}{\Gamma(t-\alpha+n/2)}
\cdot\frac{\Gamma(t)}{\Gamma(t-n)}
\\ &\qquad\times %\cdot
\biggl(\frac{\Gamma(t-\alpha)}{\Gamma(t-\alpha+n+1)}\biggr)^3
\cdot\frac{\Gamma(t-2\alpha+2n+1)}{\Gamma(t-2\alpha+n+1)}
\endalign
$$
satisfying the property
$$
t\mapsto-(t-2\alpha+n): \quad R_n(t)\mapsto(-1)^nR_n(t).
\tag1
$$
Decompose the function $R_n(t)$ as the sum of partial fractions,
$$
R_n(t)
=\sum_{j=0}^2\sum_{k=0}^n\frac{A_{jk}(n)}{(t-\alpha+k)^{3-j}},
$$
and use the last representation as in~\cite{19}, proof of Lemma~1,
to sum the quantity
$$
\align
r_n
&=\sum_{t=1}^\infty(-1)^{t-1}R_n(t)
=\sum_{t=n+1}^\infty(-1)^{t-1}R_n(t)
\tag2
\\
&=u_{0n}\sum_{\nu=1}^\infty\frac{(-1)^{n-1}}{(\nu-\alpha)^3}
+u_{1n}\sum_{\nu=1}^\infty\frac{(-1)^{n-1}}{(\nu-\alpha)^2}
+u_{2n}\sum_{\nu=1}^\infty\frac{(-1)^{n-1}}{\nu-\alpha}
-v_n,
\endalign
$$
where
$$
\gather
u_{jn}=\sum_{k=0}^n(-1)^kA_{jk}(n), \quad j=0,1,2,
\qquad
v_n=\sum_{j=0}^2\sum_{k=0}^n(-1)^nA_{jk}(n)
\sum_{\nu=1}^k\frac{(-1)^{\nu-1}}{(\nu-\alpha)^{3-j}}.
\endgather
$$
The property \thetag{1} yields $u_{0n}=u_{2n}=0$ for all $n=0,1,2,\dots$,
hence setting $u_n=u_{1n}/2$ we obtain the linear forms
$$
r_n=u_n\cdot2Z_2^-(\alpha)-v_n
\in\Bbb QZ_2^-(\alpha)+\Bbb Q,
\qquad n=0,1,2,\dots,
\tag3
$$
with effectively determined coefficients $u_n$ and $v_n$.
Applying Zeilberger's creative telescoping~\cite{12} in the
manner of~\cite{19}, Section~2 (namely, computing the
certificate and the correponding difference annihilating operator
for the function $(-1)^tR_n(t)$, which is rational
with respect to either $t$ or~$n$), we arrive at the following
recursion satisfied by both the linear forms~\thetag{3} and their
coefficients:
$$
\allowdisplaybreaks
\gather
{\align
&
(n+1)^2(n+1-\alpha)^2(5n^2-4\alpha n+\alpha^2)u_{n+1}
\\ \vspace{-1.5pt} &\;
-\bigl(55n^6-11(14\alpha-15)n^5+(179\alpha^2-385\alpha+180)n^4
\\ \vspace{-1.5pt} &\;\quad
-(116\alpha^3-358\alpha^2+332\alpha-85)n^3
+(45\alpha^4-174\alpha^3+232\alpha^2-113\alpha+15)n^2
\\ \vspace{-1.5pt} &\;\quad
-\alpha(\alpha-1)(10\alpha^3-35\alpha^2+41\alpha-12)n
+\alpha^2(\alpha-1)^2(\alpha^2-3\alpha+3)\bigr)u_n
\\ \vspace{-1.5pt} &\;
-n^2(n-\alpha)^2(5(n+1)^2-4\alpha(n+1)+\alpha^2)u_{n-1}
=0,
\tag4
\endalign}
\\
\gathered
u_0=1, \quad u_1=\alpha^2-3\alpha+3,
\\
v_0=0, \quad v_1=\frac{\alpha^2-4\alpha+5}{(\alpha-1)^2},
\endgathered
\qquad
\lim_{n\to\infty}\frac{v_n}{u_n}=2Z_2^-(\alpha).
\endgather
$$
(We justify the latter limit relation by notifying the following
consequence of the forthcoming formula~\thetag{6}:
$r_n\to0$ as $n\to\infty$.)

On the other hand, the series~\thetag{2} may be easily identified
with the very-well-poised hypergeometric ${}_6F_5(-1)$-series
$$
\alignat1
r_n
&=(-1)^nn!
\sum_{\nu=0}^\infty(3n+2-2\alpha+2\nu)
\\ &\qquad\times
\frac{\Gamma(3n+2-2\alpha+\nu)\,\Gamma(n+1-\alpha+\nu)^3\Gamma(n+1+\nu)}
{\Gamma(1+\nu)\,\Gamma(\nu+n+1-\alpha)^3\Gamma(2n+2-2\alpha+\nu)}
(-1)^\nu,
\tag5
\endalignat
$$
which admits the double integral representation
$$
r_n=u_n\cdot2Z_2^-(\alpha)-v_n
=(-1)^n\iint\limits_{[0,1]^2}
\frac{x^{n-\alpha}(1-x)^ny^{n-\alpha}(1-y)^n}{(1-x(1-y))^{n+1}}
\,\d x\,\d y
\tag6
$$
(see \cite{17}, Theorem~5).
Whipple's transformation (\cite{2}, Section~4.4, formula~(2))
gives one a more direct way to deduce the integral~\thetag{6}:
first convert the series~\thetag{5} into the hypergeometric
${}_3F_2(1)$-series
$$
r_n=\Gamma(n+1-\alpha)
\sum_{\nu=0}^\infty\frac{\Gamma(n+1+\nu)^2\Gamma(n+1-\alpha+\nu)}
{\Gamma(1+\nu)\,\Gamma(2n+2-\alpha+\nu)^2}
$$
and secondly use the Euler-type integral formula for the latter series.

\subhead
2. $Z_3(\alpha)$
\endsubhead
This time, take the rational function
$$
R_n(t)
=n!^2\cdot(2t-2\alpha+n)
\frac{(t-1)\dotsb(t-n)\cdot(t-2\alpha+n+1)\dotsb(t-2\alpha+2n)}
{\bigl((t-\alpha)(t-\alpha+1)\dotsb(t-\alpha+n)\bigr)^4}
$$
satisfying the property
$$
t\mapsto-(t-2\alpha+n): \quad R_n(t)\mapsto-R_n(t).
\tag7
$$
After partial-fraction decomposition we arrive at the quantity
$$
r_n
=\frac12\sum_{t=1}^\infty R_n(t)
=u_nZ_3(\alpha)-v_n
\in\Bbb QZ_3(\alpha)+\Bbb Q,
\qquad n=0,1,2,\dots,
$$
with effectively computable coefficients $u_n$ and $v_n$.
Zeilberger's creative telescoping produces the recursion
$$
\allowdisplaybreaks
\gather
{\align
&
(n+1)^3(n+1-\alpha)^3(2n-\alpha)(3n^2-3\alpha n+\alpha^2)u_{n+1}
\\ \vspace{-1.5pt} &\;
-(2n+1-\alpha)\bigl(102n^8-408(\alpha-1)n^7
+2(359\alpha^2-714\alpha+321)n^6
\\ \vspace{-1.5pt} &\;\quad
-6(\alpha-1)(121\alpha^2-238\alpha+83)n^5
+3(152\alpha^4-605\alpha^3+811\alpha^2-415\alpha+64)n^4
\\ \vspace{-1.5pt} &\;\quad
-2(\alpha-1)(89\alpha^4-367\alpha^3+461\alpha^2-177\alpha+15)n^3
\\ \vspace{-1.5pt} &\;\quad
+\alpha(40\alpha^5-267\alpha^4+634\alpha^3-669\alpha^2+304\alpha-45)n^2
\\ \vspace{-1.5pt} &\;\quad
-\alpha^2(\alpha-1)(2\alpha-1)(2\alpha^3-17\alpha^2+37\alpha-25)n
-\alpha^3(\alpha-1)^2(2\alpha^2-6\alpha+5)\bigr)u_n
\\ \vspace{-1.5pt} &\;
+n^3(n-\alpha)^3(2(n+1)-\alpha)(3(n+1)^2-3\alpha(n+1)+\alpha^2)u_{n-1}
=0,
\tag8
\endalign}
\\
\gathered
u_0=1, \quad u_1=2\alpha^2-6\alpha+5,
\\
v_0=0, \quad v_1=\frac{(\alpha^2-3\alpha+3)(\alpha-2)}{(\alpha-1)^3},
\endgathered
\qquad
\lim_{n\to\infty}\frac{v_n}{u_n}=Z_3(\alpha),
\endgather
$$
while writing $r_n$ as a very-well-poised hypergeometric
${}_8F_7(1)$-series and applying \cite{17}, Theorem~5,
we obtain the triple integral
$$
r_n=u_nZ_3(\alpha)-v_n
=\frac12\iiint\limits_{[0,1]^3}
\frac{x^{n-\alpha}(1-x)^ny^{n-\alpha}(1-y)^nz^{n-\alpha}(1-z)^n}
{(1-x(1-y(1-z)))^{n+1}}
\,\d x\,\d y\,\d z.
\tag9
$$
Bailey's transformation (\cite{2}, Section 6.3, formula~(2))
allows us to write the quantity~$r_n$
as the Barnes-type integral
$$
\align
r_n
&=\frac1{4\pi i}\int_{c-i\infty}^{c+i\infty}
\frac{\Gamma(n+1+s)^2\Gamma(n+1-\alpha+s)^2\Gamma(-s)^2}
{\Gamma(2n+2-\alpha+s)^2}\,\d s
\\
&=-\frac12\sum_{t=1}^\infty\frac{\d}{\d t}
\biggl(\frac{(t-1)(t-2)\dotsb(t-n)}
{(t-\alpha)(t-\alpha+1)\dotsb(t-\alpha+n)}\biggr)^2,
\endalign
$$
where the real constant $c$ lies in the interval $\Re(\alpha)-1<c<0$.
This gives another way to deduce the recursion~\thetag{8},
by applying Zeilberger's creative telescoping directly
to the latter summation (cf.~\cite{18}, Lemmas 1--3,
for the proof in the particular case $\alpha=0$). 

\subhead
3. $Z_4(\alpha)$
\endsubhead
Finally, take the rational function
$$
R_n(t)
=(2t-2\alpha+n)
\biggl(\frac{(t-1)\dotsb(t-n)\cdot(t-2\alpha+n+1)\dotsb(t-2\alpha+2n)}
{\bigl((t-\alpha)(t-\alpha+1)\dotsb(t-\alpha+n)\bigr)^2}\biggr)^2
$$
satisfying the property~\thetag{7} and consider the quantity
$$
r_n
=-\frac{(-1)^n}2\sum_{t=1}^\infty\frac{\d R_n(t)}{\d t}
=u_n\cdot6Z_4(\alpha)-v_n
\in\Bbb QZ_4(\alpha)+\Bbb Q,
\qquad n=0,1,2,\dots\,.
$$
Then Zeilberger's creative telescoping gives the recursion
$$
\allowdisplaybreaks
\gather
{\align
&
(n+1)^5(n+1-\alpha)^3(n+1-2\alpha)
(39n^4-65\alpha n^3+45\alpha^2n^2-15\alpha^3n+2\alpha^4)u_{n+1}
\\ \vspace{-1.5pt} &\;
-\bigl(10530n^{13}-1755(40\alpha-39)n^{12}
+18(11881\alpha^2-23400\alpha+10881)n^{11}
\\ \vspace{-1.5pt} &\;\quad
-9(43964\alpha^3-130691\alpha^2+122499\alpha-36036)n^{10}
\\ \vspace{-1.5pt} &\;\quad
+(497482\alpha^4-1978380\alpha^3+2795153\alpha^2-1651455\alpha+343161)n^9
\\ \vspace{-1.5pt} &\;\quad
-(449452\alpha^5-2238669\alpha^4+4229444\alpha^3-3756546\alpha^2
+1559025\alpha-241137)n^8
\\ \vspace{-1.5pt} &\;\quad
+2(149999\alpha^6-898904\alpha^5+2128142\alpha^4-2523748\alpha^3
\\ \vspace{-1.5pt} &\;\qquad
+1567577\alpha^2-480285\alpha+56394)n^7
\\ \vspace{-1.5pt} &\;\quad
-(149336\alpha^7-1049993\alpha^6+2995163\alpha^5-4449872\alpha^4
\\ \vspace{-1.5pt} &\;\qquad
+3679649\alpha^3-1676024\alpha^2+385125\alpha-33930)n^6
\\ \vspace{-1.5pt} &\;\quad
+(55088\alpha^8-448008\alpha^7+1503025\alpha^6-2693161\alpha^5+2786514\alpha^4
\\ \vspace{-1.5pt} &\;\qquad
-1681907\alpha^3+568968\alpha^2-96291\alpha+5967)n^5
\\ \vspace{-1.5pt} &\;\quad
-(14696\alpha^9-137720\alpha^8+536294\alpha^7-1132580\alpha^6+1413762\alpha^5
\\ \vspace{-1.5pt} &\;\qquad
-1065166\alpha^4+474344\alpha^3
-116539\alpha^2+13455\alpha-468)n^4
\\ \vspace{-1.5pt} &\;\quad
+\alpha(\alpha-1)(2692\alpha^8-26700\alpha^7
+105832\alpha^6-220076\alpha^5+260191\alpha^4
\\ \vspace{-1.5pt} &\;\qquad
-176174\alpha^3+65540\alpha^2-11955\alpha+780)n^3
\\ \vspace{-1.5pt} &\;\quad
-\alpha^2(\alpha-1)^2(304\alpha^7-3430\alpha^6+14198\alpha^5-29252\alpha^4
\\ \vspace{-1.5pt} &\;\qquad
+32370\alpha^3-18825\alpha^2+5265\alpha-540)n^2
\\ \vspace{-1.5pt} &\;\quad
+2\alpha^3(\alpha-1)^3(8\alpha^6-128\alpha^5+581\alpha^4-1198\alpha^3
+1220\alpha^2-558\alpha+90)n
\\ \vspace{-1.5pt} &\;\quad
+4\alpha^4(\alpha-1)^4(\alpha-2)(2\alpha-1)(\alpha^2-3\alpha+3)\bigr)u_n
\\ \vspace{-1.5pt} &\;
-n^3(n-\alpha)^3(3n-2\alpha)(3n+1-2\alpha)(3n-1-2\alpha)
\bigl(39n^4-13(5\alpha-12)n^3
\\ \vspace{-1.5pt} &\;\quad
+3(15\alpha^2-65\alpha+78)n^2-3(5\alpha^3-30\alpha^2+65\alpha-52)n
\\ \vspace{-1.5pt} &\;\quad
+(2\alpha^4-15\alpha^3+45\alpha^2-65\alpha+39)\bigr)u_{n-1}
=0,
\tag10
\endalign}
\\
\gathered
u_0=1, \quad u_1=(\alpha-1)(\alpha-2)(\alpha^2-3\alpha+3),
\\
v_0=0,
\quad v_1=-\frac{2\alpha^4-15\alpha^3+45\alpha^2-65\alpha+39}{(\alpha-1)^3},
\endgathered
\qquad
\lim_{n\to\infty}\frac{v_n}{u_n}=6Z_4(\alpha),
\endgather
$$
and Theorem~2 in~\cite{22} yields
the following $5$-fold integral:
$$
\align
r_n
&=u_n\cdot6Z_4(\alpha)-v_n
=\frac{(-1)^n\Gamma(3n+2-2\alpha)}{2\Gamma(n+1)\,\Gamma(n+1-\alpha)^2}
\\ &\qquad\times
\idotsint\limits_{[0,1]^5}
\frac{x_1^n(1-x_1)^n\prod_{j=2}^5x_j^{n-\alpha}(1-x_j)^n
\,\d x_1\dotsb\d x_5}
{(x_1(1-(1-(1-(1-x_2)x_3)x_4)x_5)+(1-x_1)x_2x_3x_4x_5)^{n+1}}.
\tag11
\endalign
$$

\subhead
4. Other recursions
\endsubhead
Polynomial recursions for approximations to the numbers
$$
Z_1^-(\alpha)=\sum_{\nu=1}^\infty\frac{(-1)^{n-1}}{\nu-\alpha},
\quad
Z_2(\alpha)=\sum_{\nu=1}^\infty\frac1{(\nu-\alpha)^2}
$$
may be constructed by means of simpler (not well-poised)
hypergeometric series. Namely, taking
$$
\align
r_n=u_nZ_1^-(\alpha)-v_n
&=(-1)^n\sum_{t=1}^\infty(-1)^{t-1}
\frac{(t-1)(t-2)\dotsb(t-n)}
{(t-\alpha)(t-\alpha+1)\dotsb(t-\alpha+n)}
\\
&=\int_0^1\frac{x^{n-\alpha}(1-x)^n}{(1+x)^{n+1}}\,\d x
\endalign
$$
we obtain the second-order recursion
$$
\allowdisplaybreaks
\gather
{\align
&
(n+1)(n+1-\alpha)(2n-\alpha)u_{n+1}
-(2n+1-\alpha)(6n^2-6(\alpha-1)n+\alpha(2\alpha-3))u_n
\\ \vspace{-1.5pt} &\qquad
+n(n-\alpha)(2(n+1)-\alpha)u_{n-1}
=0,
\tag12
\endalign}
\\
u_0=1, \quad u_1=-2\alpha+3,
\qquad v_0=0, \quad v_1=\frac{\alpha-2}{\alpha-1},
\qquad
\lim_{n\to\infty}\frac{v_n}{u_n}=Z_1^-(\alpha),
\endgather
$$
while choosing
$$
\align
r_n=u_nZ_2(\alpha)-v_n
&=(-1)^n\sum_{t=1}^\infty
\frac{n!\cdot(t-1)(t-2)\dotsb(t-n)}
{\bigl((t-\alpha)(t-\alpha+1)\dotsb(t-\alpha+n)\bigr)^2}
\\
&=(-1)^n\iint\limits_{[0,1]^2}
\frac{x^{n-\alpha}(1-x)^ny^{n-\alpha}(1-y)^n}{(1-xy)^{n+1}}
\,\d x\,\d y
\endalign
$$
we arrive at the recursion
$$
\allowdisplaybreaks
\gather
{\align
&
(n+1)^2(n+1-\alpha)^2(5n^2-6\alpha n+2\alpha^2)u_{n+1}
\\ \vspace{-1.5pt} &\;
-\bigl(55n^6-11(16\alpha-15)n^5
+2(117\alpha^2-220\alpha+90)n^4
\\ \vspace{-1.5pt} &\;\quad
-(160\alpha^3-468\alpha^2+388\alpha-85)n^3
+(56\alpha^4-240\alpha^3+316\alpha^2-142\alpha+15)n^2
\\ \vspace{-1.5pt} &\;\quad
-2\alpha(\alpha-1)(4\alpha^3-24\alpha^2+32\alpha-9)n
-2\alpha^2(\alpha-1)^2(2\alpha-3)\bigr)u_n
\\ \vspace{-1.5pt} &\;
-n^2(n-\alpha)^2(5(n+1)^2-6\alpha(n+1)+2\alpha^2)u_{n-1}
=0,
\tag13
\endalign}
\\
u_0=1, \quad u_1=-2\alpha+3,
\qquad v_0=0, \quad v_1=\frac{2\alpha^2-6\alpha+5}{(\alpha-1)^2},
\qquad
\lim_{n\to\infty}\frac{v_n}{u_n}=Z_2(\alpha).
\endgather
$$

The approach presented above allows to derive a higher-order polynomial
recursions for simultaneous approximations to odd and even zeta
values and their $\alpha$-shifts (see \cite{21}).

\subhead
5. Modular remarks
\endsubhead
It is worth mentioning that Ap\'ery's recursions for rational approximations
to $\zeta(2)$ and $\zeta(3)$ (that correspond to the case
$\alpha=0$ in~\thetag{4} or~\thetag{13} and in~\thetag{8})
have a very nice modular
interpretation: the generating function
$U(z)=\sum_{n=0}^\infty u_nz^n\in\Bbb Z[[z]]$
becomes a modular form after substituting a suitable
modular function $z=z(\tau)$ (see \cite{4, 7}).
This phenomenon happens for several other Ap\'ery-like
recursions as well (see~\cite{5, 16}); Zagier's
technique in~\cite{16} allowed to construct new simple
recursions for rational approximations to the numbers
$$
\sum_{\nu=1}^\infty\frac{\bigl(\frac{-3}\nu\bigr)}{\nu^2}
\quad\text{and}\quad
\sum_{\nu=1}^\infty\frac{\bigl(\frac{-4}\nu\bigr)}{\nu^2}
=\frac14Z_2^-\Bigl(\frac12\Bigr)
\; \text{(Catalan's constant)},
$$
where
$$
\bigl(\tfrac{-3}\nu\bigr)\equiv\nu\pmod3
\qquad\text{and}\qquad
\bigl(\tfrac{-4}\nu\bigr)\equiv\cases
0 & \text{for $\nu$ even}, \\
\nu\pmod4 & \text{for $\nu$ odd}
\endcases
$$
are quadratic characters.

In spite of the complicated form of the recursions given in Sections~1--4
above, their solutions admit nice arithmetic properties if
$\alpha$~is a rational number.
For instance, the corresponding generating functions
$U(z)=\sum_{n=0}^\infty u_nz^n$
satisfy the property $u(Az)\in\Bbb Z[[z]]$, where
the integer~$A$ depends on $\alpha\in\Bbb Q$ (although a proof of
the property in full generality for the recursions~\thetag{4}
and~\thetag{10} is still beyond reach; see \cite{10, 14, 20}
for particular results).
Beukers' computations~\cite{6} show that it is
hard to expect modular parametrizations except in the above mentioned
cases of the recursions for rational approximations to
$\zeta(2)$ and~$\zeta(3)$:
the (Zariski closures of the) Galois groups associated to the linear
differential operators annihilating functions~$U(z)$
are richer than $\SL_2$.
Beukers considers the recursion for rational approximations to~$\zeta(4)$
(the case $\alpha=0$ in~\thetag{10}) and shows that
the corresponding differential Galois group turns out to be~$\OO_5$,
while the linear differential operator corresponding to the recursion
for rational approximations
to Catalan's constant (the case $\alpha=1/2$ in~\thetag{4})
is reducible. We expect that this differential reducibility
holds for all $\alpha\notin\Bbb Z$ and that irreducible
components of the corresponding differential operators
are pullbacks of hypergeometric differential operators.
The latter fact is closely related to the general
conjecture (due to Dwork, Bombieri, $\dots$) on the structure
of differential $G$-operators. While chances to be able to attack
this general conjecture seem to be small at the moment,
to us it appears to be a nice and quite realistic program to
give a direct proof in the case of the above recursions.

\subsubhead
Acknowledgements
\endsubsubhead
This work was done during a long-term visit
at the Mathematical Institute of Cologne University.
I thank the staff of the institute and personally
P.~Bundschuh for the hospitality
and the warm working atmosphere. Special gratitude is
due to the anonymous referee of the Journal
of Computational and Applied Mathematics
for the remarks and suggestions.

%==================================================


\Refs

\ref\no1
\by R.~Ap\'ery
\paper Irrationalit\'e de $\zeta(2)$ et $\zeta(3)$
\jour Ast\'erisque
\vol61
\yr1979
\pages11--13
\endref

\ref\no2
\by W.\,N.~Bailey
\book Generalized hypergeometric series
\bookinfo Cambridge Math. Tracts
\vol32
\publ Cambridge Univ. Press
\publaddr Cambridge
\yr1935
\moreref
\bookinfo 2nd reprinted edition
\publaddr New York--London
\publ Stechert-Hafner
\yr1964
\endref

\ref\no3
\by F.~Beukers
\paper A note on the irrationality of~$\zeta(2)$ and~$\zeta(3)$
\jour Bull. London Math. Soc.
\vol11
\issue3
\yr1979
\pages268--272
\endref

\ref\no4
\by F.~Beukers
\paper Irrationality proofs using modular forms
\paperinfo Journ\'ees arithm\'etiques (Besan\c con, 1985)
\jour Ast\'erisque
\vol147--148
\yr1987
\pages271--283
\endref

\ref\no5
\by F.~Beukers
\paper On Dwork's accessory parameter problem
\jour Math. Z.
\vol241
\issue2
\yr2002
\pages425--444
\endref

\ref\no6
\by F.~Beukers
\paper Some Galois theory on Zudilin's recursions
\inbook The talk at the meeting on Elementary and Analytic Number Theory
(Mathematische Forschungsinstitut Oberwolfach, Germany,
March 9--15, 2003)
\endref

\ref\no7
\by F.~Beukers and C.\,A.\,M.~Peters
\paper A family of $K3$ surfaces and $\zeta(3)$
\jour J. Reine Angew. Math.
\vol351
\yr1984
\pages42--54
\endref

\ref\no8
\by H.~Cohen
\paper Acc\'el\'eration de la convergence
de certaines r\'ecurrences lin\'eaires
\inbook S\'eminaire de Th\'eorie des Nombres de Bordeaux (Ann\'ee 1980--81),
expos\'e 16, 2~pages
\endref

\ref\no9
\by L.\,A.~Gutnik
\paper On the irrationality of certain quantities involving~$\zeta(3)$
\jour Uspekhi Mat. Nauk [Russian Math. Surveys]
\vol34
\issue3
\yr1979
\page190
\moreref
\jour Acta Arith.
\yr1983
\vol42
\issue3
\pages255--264
\endref

\ref\no10
\by C.~Krattenthaler and T.~Rivoal
\paper Hyperg\'eom\'etrie et fonction z\^eta de Riemann
\inbook {\tt math.NT/\allowlinebreak0311114} (November 2003)
\endref

\ref\no11
\by Yu.\,V.~Nesterenko
\paper A few remarks on~$\zeta(3)$
\jour Mat. Zametki [Math. Notes]
\vol59
\yr1996
\issue6
\pages865--880
\endref

\ref\no12
\by M.~Petkov\v sek, H.\,S.~Wilf and D.~Zeilberger
\book $A=B$
\publaddr Wellesley, M.A.
\publ A.\,K.~Peters, Ltd.
\yr1996
\endref

\ref\no13
\by T.~Rivoal
\paper La fonction z\^eta de Riemann prend une infinit\'e
de valeurs irrationnelles aux entiers impairs
\jour C.~R. Acad. Sci. Paris S\'er.~I Math.
\vol331
\yr2000
\issue4
\pages267--270
\endref

\ref\no14
\by T.~Rivoal
\paper Nombres d'Euler, approximants de Pad\'e et constante de Catalan
\jour Ramanujan J.
\toappear
\endref

\ref\no15
\by V.\,N.~Sorokin
\paper One algorithm for fast calculation of~$\pi^4$
\inbook Preprint (April 2002)
\publaddr Moscow
\publ Russian Academy of Sciences,
M.\,V.~Kel\-dysh Institute for Applied Mathematics
\yr2002
\moreref
\inbook available at
{\tt http://www.wis.kuleuven.ac.be/applied/intas/Art5.pdf}
\endref

\ref\no16
\by D.~Zagier
\book Integral solutions of Ap\'ery-like recurrence equations
\bookinfo Manuscript
%\publ
%\publaddr
\yr2003
\endref

\ref\no17
\by W.~Zudilin
\paper Well-poised hypergeometric service
for diophantine problems of zeta values
\paperinfo Actes des 12\`emes rencontres arithm\'e\-tiques de Caen
(June 29--30, 2001)
\jour J. Th\'eorie Nombres Bordeaux
\vol15
\yr2003
\issue2
%\pages
\toappear
\endref

\ref\no18
\by W.~Zudilin
\paper An elementary proof of Ap\'ery's theorem
\inbook E-print {\tt math.NT/0202159} (February 2002)
\endref

\ref\no19
\by W.~Zudilin
\paper An Ap\'ery-like difference equation for Catalan's constant
\jour Electron. J. Combin.
\vol10
\issue1
\yr2003
\pages\nofrills\#R14.
\endref

\ref\no20
\by W.~Zudilin
\paper A few remarks on linear forms involving Catalan's constant
\jour Chebyshevski\v\i\ Sbornik (Tula State Pedagogical University)
\vol3
\issue2\,(4)
\yr2002
\pages60--70
\transl English transl.
\inbook E-print {\tt math.NT/\allowlinebreak0210423} (October 2002)
\endref

\ref\no21
\by W.~Zudilin
\paper A third-order Ap\'ery-like recursion for $\zeta(5)$
\jour Mat. Zametki [Math. Notes]
\vol72
\issue5
\yr2002
\pages733--737
\endref

\ref\no22
\by W.~Zudilin
\paper Well-poised hypergeometric transformations
of Euler-type multiple integrals
\inbook Pre\-print (April 2003)
\finalinfo submitted for publication
\endref

\endRefs
\enddocument